\documentclass[a4paper,12pt]{amsart}
\usepackage{amsmath,amssymb,amsthm}
\usepackage[left=2cm,right=2cm]{geometry}
\makeatletter
\@namedef{subjclassname@2020}{%
	\textup{2020} Mathematics Subject Classification}
\makeatother

\newtheorem{thm}{Theorem}
\newtheorem{res}{Result}

\newtheorem{lem}[thm]{Lemma}

\theoremstyle{definition}
\newtheorem{defn}{Definition}
\theoremstyle{remark}

\def\C{\mathbb C}

\def\hol{\mathcal O}

\def\O{\Omega}

\def\eps{\varepsilon}

\def\k{\kappa}
\def\endofproof{\hfill \square}

\def\lk{ l^{\kappa}}
\def\Cn{\mathbb{C}^n}

\renewcommand{\Im}{\operatorname{Im}}
\renewcommand{\Re}{\operatorname{Re}}

\begin{document}
	
	\title[Visible $\mathcal C^2$-smooth domains are pseudoconvex]
	{Visible $\mathcal C^2$-smooth domains are pseudoconvex}

	\author{Nikolai Nikolov, Ahmed Yekta \"Okten, Pascal J. Thomas} 
	
	\address{N. Nikolov\\
		Institute of Mathematics and Informatics\\
		Bulgarian Academy of Sciences\\
		Acad. G. Bonchev Str., Block 8\\
		1113 Sofia, Bulgaria}
	
	\address{Faculty of Information Sciences\\
		State University of Library Studies
		and Information Technologies\\
		69A, Shipchenski prohod Str.\\
		1574 Sofia, Bulgaria}
	
	\email{nik@math.bas.bg}
	
	\address{A. Y. \"Okten\\
		Institut de Math\'ematiques de Toulouse; UMR5219 \\
		Universit\'e de Toulouse; CNRS \\
		UPS, F-31062 Toulouse Cedex 9, France} \email{ahmed$\_$yekta.okten@math.univ-toulouse.fr}

	\address{P. J. Thomas\\	Institut de Math\'ematiques de Toulouse; UMR5219 \\Universit\'e de Toulouse; CNRS \\	UPS, F-31062 Toulouse Cedex 9, France} \email{pascal.thomas@math.univ-toulouse.fr}
	\thanks{The first author was partially supported by the Bulgarian National Science Fund, Ministry of Education
		and Science of Bulgaria under contract KP-06-N52/3. The second author received support from the University
		Research School EUR-MINT (State support managed by the National Research Agency for Future Investments
		program bearing the reference ANR-18-EURE-0023). The third author acknowledges 
		the hospitality of the Vietnam Institute for Advanced Studies in Mathematics
		for part of the time during which this work was carried out.}

	\begin{abstract}  We show that a domain in $\mathbb C^n$ with $\mathcal C^2$-smooth boundary which satisfies the visibility property is pseudoconvex.  
	\end{abstract}
	
	
	
	\maketitle
\section{Introduction and statement of the result}

Various holomorphic invariants are used to understand the properties of domains 
in $\mathbb C^n$ (or indeed complex manifolds or spaces) and their mappings. Among them are
those infinitesimal Finsler metrics which are decreasing under
holomorphic maps, and the distances obtained from them. The
largest and best known of those is the Kobayashi metric. 

In the last couple of decades, interest has grown in the study of the metric geometric properties 
of domains in $\mathbb C^n$ when endowed with the Kobayashi metric. Visibility is
a property of the domain as a metric space, and of its boundary---under a specific embedding
in the Euclidean space $\mathbb C^n$. Visibility will be defined 
precisely below. Intuitively it means that 
near-geodesics  (curves that almost minimize length) between two points 
close to two distinct boundary points
have to ``curve back'' and meet some relatively compact subset of the domain
depending only on the two boundary points.  
The visibility property is well-studied and has many applications, see for instance \cite{BZ}, \cite{BNT}, \cite{BM} and \cite{S}.

Visibility  clearly holds when the domain is the unit ball (and does not hold for the polydisc);
general considerations about Gromov hyperbolic metric spaces show that it holds 
for Gromov hyperbolic domains when their Euclidean
boundary can be identified with the Gromov boundary, which is the case for $\mathcal C^2$-smooth 
strongly pseudoconvex domains \cite{BB} or smooth convex domains of finite type \cite{Z}. 

Early examples of domains that satisfy the visibility property, in particular the Goldilocks domains introduced by Bharali and Zimmer \cite[Definition 1.1]{BZ}, were known to be pseudoconvex. Recently, in \cite{Ban}, the author exhibited non-pseudoconvex domains that satisfy the visibility property. However, the examples provided in \cite[Theorems 1.2, 1.4 and 3.1]{Ban} have quite irregular boundaries:
they are disconnected and contain portions which are ``small'', precisely of real codimension 
at least $2$. In the same paper, the author mentions a question of Filippo Bracci 
\cite[Question 5.1]{Ban}, asking
whether a similar example with a $\mathcal C^1$ boundary could be found. 
We answer this in the negative, but only in the more restrictive $\mathcal C^2$ case:
a domain with $\mathcal C^2$-smooth boundary which enjoys
the visibility property must be pseudoconvex.

To be more precise, let $\Omega$ be a domain in $\C^n$, $z,w\in \Omega$ and $v\in\Cn$. Recall that the \emph{Kobayashi pseudodistance $k_\Omega$} is the largest pseudodistance which does not exceed the Lempert function	
\[
l_\O(z,w):=\inf\{\tanh^{-1}|\alpha|:\exists\varphi\in\hol(\Delta,\O)	\hbox{ with }\varphi(0)=z,\varphi(\alpha)=w\},
\]
where $\Delta$ is the unit disc and $\mathcal{O}(M,N)$ denotes the space of holomorphic functions defined on a complex manifold $M$ into a complex manifold $N$.
	
Also recall the definition of \emph{Kobayashi-Royden pseudometric,}
	$$
	\k_\O(z;v)=\inf\{|\alpha|:\exists\varphi\in\mathcal{O}(\Delta,\O)\hbox{ with }
	\varphi(0)=z,\alpha\varphi'(0)=v\}.
	$$

The \emph{Kobayashi-Royden length} of an absolutely continuous curve $\gamma:[0,l]\rightarrow \O$ is defined as 
$$ 
\lk_\O(\gamma):=\int_0^l \k_\O(\gamma(t);\gamma'(t))dt. 
$$
By \cite[Theorem 1.2]{V}, it turns out that $k_\Omega$ is the inner distance associated with the Kobayashi-Royden pseudometric. That is 
\begin{equation}
\label{inflength}
k_\O(z,w):=\inf \left\{ \lk_\O(\gamma) \mbox{ where }\gamma
\mbox{ absolutely continuous curve joining }z\mbox{ to }w \right\}.
\end{equation}
	
For $\lambda \ge 1$, $\eps\ge 0$, an absolutely continuous curve $\gamma: [0,l]\rightarrow\O$ is said to be a 
\emph{$(\lambda,\epsilon)$-geodesic} if for all $ t_1,t_2 \in [0,l]$ we have that 
$$ 
\lk_\O(\gamma|_{[t_1,t_2]}) \leq \lambda k_\O(\gamma(t_1),\gamma(t_2))+\epsilon. 
$$
In this terminology, geodesics with respect to the Kobayashi-Royden 
infinitesimal metrics are $(1,0)$-geodesics, equivalently, they 
attain the infimum in \eqref{inflength}. The existence of such a minimum
is not always guaranteed,  
while the definition of the Kobayashi distance as an infimum ensures 
the existence of $(1,\epsilon)$-geodesics for any $\epsilon>0$. So it is useful to consider the wider notions above, introduced in \cite{BZ}. 

\begin{defn}
A domain $\Omega\subset\Cn$ satisfies the \emph{$(\lambda,\epsilon)$-visibility property} if for any pair of distinct points $p,q \in \partial \Omega$  there exist neighborhoods $U,V$ of $p,q$ respectively such that $\overline{U}\cap \overline{V} = \emptyset$, and a compact set $K:=K_{p,q,\lambda,\epsilon}\subset\subset \Omega$ such that if $\gamma:[0,l]\to \Omega$ is a $(\lambda,\epsilon)$-geodesic with $\gamma(0)\in\O\cap U$ and $\gamma(l)\in\O\cap V$ then $\gamma([0,l])\cap K \neq \emptyset$. 

The domain $\Omega$ satisfies \emph{the visibility property} if it satisfies the $(\lambda,\epsilon)$-visibility property for any $\lambda \geq 1$ and $\epsilon \geq 0$. 
\end{defn}

\begin{thm}\label{thm:profthomassresult}
Let $\Omega$ be a bounded domain in $\mathbb C^n$ with $\mathcal C^2$-smooth boundary. 
Suppose that either
\begin{enumerate}
\item
$\Omega$ satisfies the $(1,\epsilon)$-visibility property for some $\epsilon>0$; or 
\item
$\Omega$ satisfies the $(\lambda,0)$-visibility property for all $\lambda > 1$.
\end{enumerate}
	Then $\Omega$ is pseudoconvex. \end{thm}   


\section{Proof of Theorem \ref{thm:profthomassresult}}

Let $\O$ be a domain in $\Cn$ with $\mathcal C^2$-smooth boundary. Set $\delta_\O(z):=\min_{w\in\partial\O}\|z-w\|$ to be the boundary distance function and let $d_\O$ denote the signed boundary distance function, that is 
$d_\O(z):= -\delta_\O(z)$ if $z\in\O$ and $d_\O(z):= \delta_\O(z)$ otherwise.

\begin{lem}\label{lem:signeddistanceandprojection}\cite[Lemma 2.1]{BB}
	Suppose that $\Omega$ is a bounded domain in $\Cn$ with $\mathcal{C}^{2}$-smooth boundary. Let $(a,b)$ denote the line segment joining $a,b\in\Cn$ and for $\eta>0$ set $N_\eta:= \cup_{p\in\partial\O}(p-\eta\nu_p,p+\eta\nu_p)$,
where $\nu_p$ is the inner unit normal to $\partial\Omega$ taken at the point $p$. Then there exists a small enough $\eta>0$ such that:	
	
	\emph{(i)} For all $z\in N_\eta$ there exists a unique point $\pi_\O(z)\in\partial\O$ such that $\|z-\pi_\O(z)\|=\delta_\O(z)$.
	
	\emph{(ii)} $d_\O:\Cn \to \mathbb{R}$ is $\mathcal{C}^2$-smooth on $N_\eta$. 
	
	\emph{(iii)} For $z\in N_\eta$ the signed boundary distance function satisfies $2\overline{\partial} d_\O|_z= 2\overline{\partial} d_\O|_{\pi_\O(z)}=-\nu_z$,
	where we write $\nu_z:= \nu_{\pi_\O(z)}$.
	
	\emph{(iv)} $\pi_\O:\Cn \to \mathbb{R}$ is $\mathcal{C}^1$-smooth on $N_\eta$ and for any $p\in\partial\O$ the fibers of this map satisfy $\pi^{-1}_\O(p) \supset (p-\eta\nu_p,p+\eta\nu_p)$.
\end{lem}

Let $\O$ be a domain in $\mathbb C^n$ with $\mathcal C^2$-smooth boundary, $z\in N_\eta$ and $v \in \mathbb C^n$. Denote the standard Hermitian inner product by $\langle z,w \rangle_{\mathbb C}:=
\sum_{j=1}^n z_j\bar w_j$. 
At the basepoint $z\in N_\eta$, we write a vector $v$ in $\mathbb C^n$ 
as $v=v_H+v_N$ where $v_N:=  \langle v, \nu_z \rangle_{\mathbb C} \nu_z$ and $v_H:=v-v_N$. 
The component $v_H$ is known as complex-tangential or \emph{horizontal}. 

The following estimates relate the behavior of a $\mathcal C^1$-smooth curve and of its projection to the boundary.

\begin{lem}\cite[Lemma 2.2]{BB}\label{lem:comparinginfinitesimallengths}
	Let $\gamma:[0,l]\to N_\eta$ be a $\mathcal{C}^1$-smooth curve and $\alpha:=\pi_\O\circ \gamma$. Then there exists a constant $C>0$ such that the following estimates hold: 
	
	\emph{(i)} $\|(\gamma'(t))_H-(\alpha'(t))_H\| \leq C \delta_\O(\gamma(t))\|\alpha'(t)\|$. 
	
	\emph{(ii)} $\|(\gamma'(t))_N\| \leq  \|(\alpha'(t))_N\| + C\delta_0\|\alpha'(t)\|$ if in addition $\delta_\O(\gamma(t))=\delta_0$ for all $t\in [0,l]$.
\end{lem}

Let $\Omega$ be a bounded domain in $\Cn$ with $\mathcal C^2$-smooth boundary. Recall that $p\in\partial \Omega$ is a \emph{non-pseudoconvex boundary point} if the restriction of the Levi form of $\Omega$ at the point $p$ has at least one negative eigenvalue. 
If $\O\subset\mathbb C^2$, observe that $p\in\partial\Omega$ is a non-pseudoconvex boundary point if and only if  $\mathbb C^2\setminus\overline \Omega$ is strongly pseudoconvex near $p$. The growth of the Kobayashi-Royden pseudometric near non-pseudoconvex boundary points has been studied in \cite{DNT}.

\begin{lem}\cite[Proposition 3]{DNT}\label{lem:DNTlemma}
	Let $\O$ be a domain in $\mathbb{C}^2$ with $\mathcal C^2$-smooth boundary, and $p\in\partial\Omega$ be a non-pseudoconvex boundary point. Then, there exists $C_L>0$ such that
	$$ \kappa_\O(z;v) \leq C_L\left(\frac{\|v_N\|}{\delta^{3/4}_\Omega(z)}+\|v\|\right) \:\:\:\:\: \text{for} \:\:\: z\in\O \:\:\: \text{near} \:\:\: p, \:\:\: v \in \mathbb{C}^2. $$
\end{lem}


\begin{proof}[Proof of Theorem \ref{thm:profthomassresult}]

Let $\O:=\{\rho(z)<0\}\subset\Cn$ be a non-pseudoconvex domain with $\mathcal C^2$-smooth boundary, and let $p\in\partial\Omega$ be a non-pseudoconvex boundary point. For each $\epsilon>0$, we will find points $p, q \in \partial\Omega$ which fail the $(1,\epsilon)$-visibility property . 
	
	By taking affine transformations, we assume that $p=0$ and $\nu_p=(-1,0,...,0)$.	Since $p=0$ is not a Levi pseudoconvex boundary point, there exists a vector $v\in\mathcal \{0\}\times\mathbb{C}^{n-1} \setminus \{0\}$ such that the Levi form of $\rho$ at $p$ satisfies $\mathcal{L_{\rho}}(v,v) <0$. By taking a rotation in $\{0\}\times \mathbb C^{n-1}$ we also assume that $v:=(0,1,0,...,0)$. Choose a small enough neighborhood $U$ of $p$ such that $\Omega' \subset \mathbb{C}^2$ given by $\Omega':=\{(z_1,z_2):(z_1,z_2,0,...,0)\in\Omega\cap U\}$ is a domain in $\mathbb{C}^2$. Note that this is possible because $\Omega$ has a $\mathcal{C}^2$-smooth boundary, and hence its boundary is locally connected. By choosing $U$ appropriately, we may furthermore assume that $\Omega'$ has a $\mathcal C^2$-smooth boundary.
	
Observe that: 

(a) $\Omega'$ near $p':=(0,0)\in\partial \Omega'$ is given by $$\Omega':=\{\rho(z_1,z_2,0,...,0)<0\},$$ hence $p'$ is a non-pseudoconvex boundary point of $\Omega'$. 

(b) The map $i:\Omega'\to\Omega$ given by $i(z_1,z_2)=(z_1,z_2,0,...,0)$ is a holomorphic embedding.

Since $p'$ is a non-pseudoconvex boundary point of $\Omega'$,
we can choose a smaller neighborhood $U'$ of $p$
so that $U'\setminus \Omega'$ is strongly pseudoconvex.
Reducing $U'$ if needed, we may fix an $\eta_0>0$ 
such that the conclusion of Lemma \ref{lem:DNTlemma} holds on the open set $ N':= \left(\cup_{p\in\partial\O'\cap U'} (p-\eta_0\nu_{p'},p+\eta_0\nu_{p'}) \right) \cap \Omega'$, $\nu_{p'}=(-1,0)$.

To prove part (1) of the theorem we recall the following result of Chow.

\begin{res}\label{res:chow}\cite{Cho}\cite[Theorem 2.4, p. 15]{Bel}
Let $M$ be a connected Riemannian manifold and $S:=\{X_1,...,X_N\}$ be a set of $\mathcal{C}^1$-smooth vector fields on $M$. Suppose that the iterated Lie brackets of the elements of $S$ generate the (real) tangent space $T_p M$ at any $p\in\partial\Omega$. Then, any $x,y\in M$ can be joined by an integral curve $\alpha:[0,1]\to M$ of a vector field $X$, where for any $t\in[0,1]$, the vector field $X$ at $\alpha(t)$ belongs to the span of the elements of $S$. 
\end{res} 

As the $U'\setminus \Omega'$  is strongly pseudoconvex, 
the result above implies  that 
for any $q'\in\partial \O' \cap U'$ we may find a complex tangential $\mathcal C^1$-smooth curve that connects $p'$ to $q'$, that is, 
a $\mathcal C^1$-smooth curve $\alpha:[0,l]\to \partial \Omega'$ with $\alpha(0)=p'$, $\alpha(l)=q'$, $\alpha'(t)=(\alpha'(t))_H$ for any $t\in [0,l]$. 

Recall that the \emph{Carnot-Carath\'eodory distance} $d_{CC}(p,q)$ is defined as the infimum
of the lengths of the horizontal, i.e. complex tangential, rectifiable curves connecting $p$ to $q$.
Here we use Euclidean length which, as pointed out in \cite[(1.1), p. 506]{BB}, is equivalent
up to multiplicative constants to the definition using the Levi form. 
The following estimate for the Carnot-Carath\'eodory distance on strongly pseudoconvex domains is called the ball-box estimate in \cite[Proposition 3.1]{BB}.

\begin{res}\label{res:ballbox}\cite[Proposition 3.1]{BB}
	Let $D$ be a bounded strongly pseudoconvex domain with $\mathcal C^2$-smooth boundary. There exists $\epsilon_0>0$ and $C>1$ such that for all $\epsilon \in (0,\epsilon_0)$ and $p\in\partial D$ we have
	$$ 
	\text{Box}(p,\epsilon/C) \leq B_{CC}(p,\epsilon) \leq \text{Box}(p,C \epsilon), 
	$$ 
where $\text{Box}(p,r):=\{p+v \in\partial D: \|v_H\|< r, \|v_N\|<r^2\}$ 
\newline
and $B_{CC}(p,r):=\{x\in\partial D: d_{CC}(p,x) < r\}$.
\end{res}

Let  $l_e(\gamma)$ denotes the Euclidean length of a curve $\gamma$.
As a consequence of the above result applied to $U'\setminus \Omega'$ there exists $C'>0$ such that,
when  $\alpha:[0,l] \to \partial \Omega'$ is a piecewise $\mathcal C^1$ curve with Euclidean length  approximating the Carnot-Carath\'eodory distance from $p'$ to $q'$, 
\begin{equation}
\label{eqn:ballbox} 
l_e(\alpha) \leq C'(\|p'-q'\|+ |\langle p'-q', \nu_{p'} \rangle_{\mathbb C}|^{1/2}).
\end{equation}

\emph{Claim.} Let $\alpha_\eta(t):=\alpha(t)+\eta \nu_{\alpha(t)}$. There exists $K>0$ such that for small enough $\eta>0$ we have $\lk_{\O'}(\alpha_\eta) \leq K l_e(\alpha).$

\emph{Subproof of Claim.} We may assume that for small enough $\eta>0$, $\alpha_\eta(t)$ remains in $N'$, hence the claim immediately follows from Lemmas \ref{lem:comparinginfinitesimallengths} and \ref{lem:DNTlemma}. $\endofproof$

By \eqref{eqn:ballbox} and the monotonicity of the Kobayashi-Royden pseudometric under holomorphic maps, our claim gives
$$
\lk_{\Omega}(\alpha_\eta|_{[t_1,t_2]})\leq 
\lk_{\Omega}(\alpha_\eta)\leq \lk_{\Omega'}(\alpha_\eta) \leq K l_e(\alpha) \leq  
2 K C'\|p'-q'\|^{1/2} \leq k_\O(\alpha_\eta(t_1),\alpha_\eta(t_2))+ 2 K C'\|p'-q'\|^{1/2} .
$$
Our construction shows that
$\alpha_\eta$ are $(1, c(p',q'))$-geodesics, where $c(p',q'):= 2 K C'\|p'-q'\|^{1/2}$. Moreover $\max_{t\in[0,l]}\delta_\Omega(\alpha_\eta(t)) \leq \eta$. By taking $q'$ close enough to $p'$, 
$c(p',q') < \epsilon$, and letting $\eta\to 0$ the theorem follows. 

Observe that we could choose any two points close enough to $p$ to violate
 $(1,\epsilon)$-visibility.

To prove part (2) of the Theorem, we will construct special curves such 
that any arc on the curve verifies that its Kobayashi length is comparable to 
its Euclidean length, itself comparable to the Euclidean distance between its extremities.

Recall that since $\O$ is bounded, there exists a constant $C_\O >0$ 
such that for any $z\in \O$, $v\in \C^n$, $\k (z;v) \ge C_\O \|v\|$,
and thus for any rectifiable curve $\gamma:[a,b]\longrightarrow \O$, 
$\lk_\O (\gamma) \ge  C_\O l_e (\gamma) \ge
 C_\O \|\gamma(a)-\gamma(b)\|$; passing to the infimum, 
$k_\O (\gamma(a),\gamma(b)) \ge C_\O \|\gamma(a)-\gamma(b)\|$.

Choose a $\mathcal C^1$ vector field $v:\partial \O' \longrightarrow \C^2$ such that 
for any $\zeta \in \partial \O'$, $v(\zeta) \in T^{\C}_\zeta \partial \O'$ (the complex
tangent space to $\partial \O'$ at $\zeta$), $\|v(\zeta)\|=1$, and $v(p)= (0,1)$. 
This can be done by choosing at each point $\zeta$ in a small enough neighborhood
of $p$ the unique unit vector  $(v_1,v_2) \in T_\zeta \partial \O' \cap i T_\zeta \partial \O'
\cap \{ \Im v_2=0\}$ satisfying $\Re v_2 >0$; it depends $\mathcal C^1$-smoothly
on $\zeta$ because $\partial \O'$ is $\mathcal C^2$-smooth. 

Let $\alpha$ be an integral curve of $v$ verifying $\alpha(0)=p$, which we restrict to 
the interval $[0,s]$.  By construction, $\alpha$ will be $\mathcal C^2$-smooth, 
and $\alpha'(p)=v(p)=(0,1)$. Thus $\|\alpha'(t)-(0,1)\| \le C|t|$ and
for $t_1, t_2$ small enough, 
\[
\alpha(t_1)-\alpha(t_2)= (0,t_1-t_2) + O(|t_1-t_2|^2).
\]
Therefore, given any $\epsilon >0$, we can choose $s$ small enough so that for $0\le t_1 < t_2\le s$,
\[
\| \alpha(t_1) - \alpha(t_2) \| \le l_e(\alpha|_{[t_1,t_2]})
= |t_1-t_2| \le (1+\epsilon) \| \alpha(t_1) - \alpha(t_2) \| .
\]
Define $\alpha_\eta$ as in the Claim above. Then by reducing $s$ and taking $\eta$
small enough, $\alpha_\eta$ verifies 
\[
 \| \alpha_\eta(t_1) - \alpha_\eta(t_2) \| \le l_e(\alpha_\eta|_{[t_1,t_2]})
\le  (1+\epsilon) |t_1-t_2| \le (1+\epsilon)^2 \| \alpha_\eta(t_1) - \alpha_\eta(t_2) \| .
\]
By Lemmas \ref{lem:comparinginfinitesimallengths} and \ref{lem:DNTlemma}, 
$\k_{\O'} (\alpha_\eta(t) ; \alpha'_\eta(t) ) \le C_L (1+O(\eta)) (1+ O(\eta^{1/4})) \le C_1$.
Therefore
\[
\lk_{\Omega}(\alpha_\eta|_{[t_1,t_2]})\leq C_1 l_e(\alpha_\eta|_{[t_1,t_2]})
\le C_1 (1+\epsilon)^2 \| \alpha_\eta(t_1) - \alpha_\eta(t_2) \|  
\le C_\O^{-1} C_1 (1+\epsilon)^2 k_\O (\alpha_\eta(t_1) , \alpha_\eta(t_2)),
\]
so for $\lambda > C_\O^{-1} C_1 (1+\epsilon)^2$,  $(\lambda,0)$-visibility
is violated.
\end{proof}
\noindent
{\bf Remarks.} 
(a) It follows from the proof that for the family of curves in the proof of Part (2), 
\[
 \| \alpha_\eta(t_1) - \alpha_\eta(t_2) \| \asymp l_e(\alpha_\eta|_{[t_1,t_2]}) \asymp
\lk_{\Omega}(\alpha_\eta|_{[t_1,t_2]}) \asymp k_\O (\alpha_\eta(t_1) , \alpha_\eta(t_2)).
\]
\noindent
(b) In the case where $n=2$, by using the estimates \eqref{eqn:ballbox} 
one could see that, as is the case
for $(1,\epsilon)$-geodesics, for 
any $p, q$ close enough to $(0,0)$, one can find a family of curves tending to
$\partial \Omega$ connecting $p_k, q_k$
with $p_k\to p$ and $q_k\to q$, which are $(\lambda,0)$-geodesics for some
large $\lambda$ depending on $p$ and $q$.

\vskip.3cm
The authors wish to thank the anonymous referee for numerous comments that greatly improved the exposition.

	{}
	
\end{document}